\documentclass[a4paper]{amsart}
\usepackage[OT2,T1]{fontenc}
\usepackage[dvips]{graphicx}
\theoremstyle{plain}
  \newtheorem{thm}{Theorem}[section]

  \newtheorem{conj}[thm]{Conjecture}

\theoremstyle{definition}

\theoremstyle{remark}
  \newtheorem{rem}[thm]{Remark}
  \newtheorem*{ack}{Acknowledgments}
\DeclareMathAlphabet{\cyr}{OT2}{wncyr}{m}{n}
\newcommand{\Lob}{\operatorname{\cyr{L}}}

\newcommand{\Z}{\mathbb{Z}}
\newcommand{\C}{\mathbb{C}}

\newcommand{\darrows}{\swarrow\hspace{-1em}\searrow}
\newcommand{\uarrows}{\nwarrow\hspace{-1em}\nearrow}
\newcommand{\rarrows}{\nearrow\hspace{-1em}\searrow}
\newcommand{\larrows}{\swarrow\hspace{-1em}\nwarrow}
\newcommand{\qbar}{\overline{q}}
\newcommand{\qbinom}[2]{\genfrac{[}{]}{0pt}{}{#1}{#2}}

\newcommand{\Li}{\operatorname{Li}}
\newcommand{\Vol}{\operatorname{Vol}}

\numberwithin{equation}{section}
\date{\today}
\allowdisplaybreaks
\begin{document}
\date{6 April, 2000}
\title[The colored Jones function of a knot and its volume]
{The asymptotic behavior of the colored Jones function of a knot and its volume}
\author{Hitoshi Murakami}
\address{
  Department of Mathematics,
  School of Science and Engineering,
  Waseda University,
  Ohkubo, Shinjuku, Tokyo 169-8555, Japan
}
\curraddr{
Department of Mathematics,
Tokyo Institute of Technology,
Oh-okayama, Meguro, Tokyo 152-8551, Japan
}
\email{starshea@tky3.3web.ne.jp}
\begin{abstract}
I show various calculations of the limit of the colored Jones function for the
figure-eight knot and confirm R.~Kashaev's conjecture in this case.
\end{abstract}
\maketitle
\section{Introduction}
The Jones polynomial $V(L;t)\in\Z[t^{1/2},t^{-1/2}]$ was first introduced by V.~Jones \cite{Jones:BULAM3} as a link invariant which satisfies the following recursive formula:
\begin{equation*}
  \begin{cases}
    V(O;t)=1,
    \\
    tV(L_{+};t)-t^{-1}V(L_{-};t)=(t^{1/2}-t^{-1/2})V(L_{0};t),
  \end{cases}
\end{equation*}
where $O$ is the trivial knot and $(L_{+},L_{-},L_{0})$ is a usual skein triple of links.
(The original version uses a different normalization.)
Let $V_N(L;t)$ denote the colored Jones polynomial, colored with the $N$-dimensional irreducible representation of $sl_2(\C)$, which can be defined
by using an enhanced Yang-Baxter operator \cite{Turaev:INVEM88}.
We use the normalization for $V_N(L;t)$ such that $V_2(L;t)=V(L;t)$
(Note that $V_N(O;t)=1$.)
We also denote $V_N(L;\exp(2\pi\sqrt{-1}/N))$ by $J_N(L)$ and call it the colored Jones function.
\par
In \cite{Murakami/Murakami:volume} we proved that $J_N(L)$ is the same as Kashaev's invariant \cite{Kashaev:LETMP97} and generalized his conjecture to
\begin{conj}[Volume Conjecture]
  \begin{equation*}
    \frac{2\pi}{v_3}\lim_{N\to\infty}\frac{\log|J_N(K)|}{N}=\Vert{K}\Vert
  \end{equation*}
for any knot $K\in S^3$.
Here $v_3$ is the volume of the ideal regular hyperbolic tetrahedron and
$\Vert{K}\Vert$ is the Gromov norm of $S^3\setminus{K}$ \cite{Gromov:INSHE82}.
\end{conj}
\begin{rem}
Kashaev's conjecture is for hyperbolic knots;
\begin{conj}[Kashaev \cite{Kashaev:LETMP97}]
  Let $K$ be a hyperbolic knot then
  \begin{equation*}
    2\pi\lim_{N\to\infty}\frac{\log|J_N(K)|}{N}=\Vol(S^3\setminus K).
  \end{equation*}
\end{conj}
Note that since for a hyperbolic knot $\Vert{K}\Vert=\Vol(S^3\setminus K)/v_3$,
our conjecture is a generalization of Kashaev's conjecture.
\end{rem}
\par
In this article I will describe how to calculate the colored Jones function
and show various ways to confirm Kashaev's conjecture for the figure-eight knot.
\begin{rem}
\quad\par
\begin{enumerate}
\item
  Kashaev proved his conjecture for the knots $4_1$, $5_2$, and $6_1$
  \cite{Kashaev:LETMP97}.
\item
  Kashaev and O.~Tirkkonen proved the volume conjecture for torus knots
  \cite{Kashaev/Tirkkonen:1999}.
\item
  Y.~Yokota proved Kashaev's conjecture for the knot $6_2$ in
  \cite{Yokota:Murasugi70} and suggests a proof for general hyperbolic knots.
  See also his forthcoming paper \cite{Yokota:volume}.
\item
  The volume conjecture is proved for the knots $6_3$, $7_2$ and $8_9$, and for the Whitehead link in \cite{Murakami/Murakami/Okamoto/Takata/Yokota:volume}.
\end{enumerate}
\end{rem}
\begin{ack}
I would like to thank the participants in the workshops held at the International Institute for Advanced Study in October 1999 (some information is available on the WWW at http://www.iias.or.jp/research/suuken/1999.10.12.html) and at Kyushu University in December 1999.
Especially I am grateful to J.~Murakami, M.~Okamoto, T.~Takata, and Y.~Yokota for useful conversations.
\end{ack}
\section{How to calculate the colored Jones function from a $(1,1)$-tangle}
In this section I show how to calculate the colored Jones function
of a knot from its $(1,1)$-tangle diagram.
\par
Let $D$ be a $(1,1)$-tangle diagram of a given knot $K$.
We assume that $D$ is generic with respect to the height function with end
points at the top and at the bottom.
In particular each crossing in $D$ is one of the following eight forms.
The diagram is decomposed into arcs after we cut it into four arcs at each
crossing.
We label each arc and each crossing with parameter from $0$ to $N-1$.
\par
For each crossing we assign a complex number as follows, where $\delta_{i,j}$ is
Kronecker's delta.
\setlength{\unitlength}{1mm}
\medskip
\begin{align*}
\raisebox{-5mm}{\begin{picture}(10,10)
  \put(10,10){\vector(-1,-1){10}}
  \put( 6, 4){\vector( 1,-1){ 4}}
  \put( 4, 6){\line(-1, 1){ 4}}
  \put( 0,10){\makebox(0,0)[br]{$i$}}\put(10,10){\makebox(0,0)[bl]{$j$}}
  \put( 0, 0){\makebox(0,0)[tr]{$k$}}\put(10, 0){\makebox(0,0)[tl]{$l$}}
  \put( 5, 5){\makebox(0,0)[l]{\,\,${\boldsymbol m}$}}
\end{picture}}
\quad &:\,
  (\darrows_{m}^{+})_{k,l}^{i,j}:=
  \delta_{l,i+m}\delta_{k,j-m}
  \frac{(q)_{i}^{-1}(q)_{j}}
       {(q)_{m}(q)_{k}(q)_{l}^{-1}}
  (-1)^{i+k+1}
  q^{ik+(i+k)/2+(N^2+1)/4}
  \\[8mm]
\raisebox{-5mm}{\begin{picture}(10,10)
  \put( 0, 0){\vector( 1, 1){10}}
  \put( 6, 4){\line( 1,-1){ 4}}
  \put( 4, 6){\vector(-1, 1){ 4}}
  \put( 0,10){\makebox(0,0)[br]{$i$}}\put(10,10){\makebox(0,0)[bl]{$j$}}
  \put( 0, 0){\makebox(0,0)[tr]{$k$}}\put(10, 0){\makebox(0,0)[tl]{$l$}}
  \put( 5, 5){\makebox(0,0)[l]{\,\,${\boldsymbol m}$}}
\end{picture}}
\quad &:\,
  (\uarrows_{m}^{+})_{k,l}^{i,j}:=
  \delta_{i,l+m}\delta_{j,k-m}
  \frac{(q)_{i}(q)_{j}^{-1}}
       {(q)_{m}(q)_{k}^{-1}(q)_{l}}
  (-1)^{j+l+1}
  q^{jl+(j+l)/2+(N^2+1)/4}
  \\[8mm]
\raisebox{-5mm}{\begin{picture}(10,10)
  \put( 0,10){\vector( 1,-1){10}}
  \put(10,10){\line(-1,-1){ 4}}
  \put( 4, 4){\vector(-1,-1){ 4}}
  \put( 0,10){\makebox(0,0)[br]{$i$}}\put(10,10){\makebox(0,0)[bl]{$j$}}
  \put( 0, 0){\makebox(0,0)[tr]{$k$}}\put(10, 0){\makebox(0,0)[tl]{$l$}}
  \put( 5, 5){\makebox(0,0)[l]{\,\,${\boldsymbol m}$}}
\end{picture}}
\quad &:\,
  (\darrows_{m}^{-})_{k,l}^{i,j}:=
  \delta_{l,i-m}\delta_{k,j+m}
  \frac{(\qbar)_{i}(\qbar)_{j}^{-1}}
       {(\qbar)_{m}(\qbar)_{k}^{-1}(\qbar)_{l}}
  (-1)^{j+l+1}
  q^{-jl-(j+l)/2-(N^2+1)/4}
  \\[8mm]
\raisebox{-5mm}{\begin{picture}(10,10)
  \put(10, 0){\vector(-1, 1){10}}
  \put( 6, 6){\vector( 1, 1){ 4}}
  \put( 4, 4){\line(-1,-1){ 4}}
  \put( 0,10){\makebox(0,0)[br]{$i$}}\put(10,10){\makebox(0,0)[bl]{$j$}}
  \put( 0, 0){\makebox(0,0)[tr]{$k$}}\put(10, 0){\makebox(0,0)[tl]{$l$}}
  \put( 5, 5){\makebox(0,0)[l]{\,\,${\boldsymbol m}$}}
\end{picture}}
\quad &:\,
  (\uarrows_{m}^{-})_{k,l}^{i,j}:=
  \delta_{i,l-m}\delta_{j,k+m}
  \frac{(\qbar)_{i}^{-1}(\qbar)_{j}}
       {(\qbar)_{m}(\qbar)_{k}(\qbar)_{l}^{-1}}
  (-1)^{i+k+1}
  q^{-ik-(i+k)/2-(N^2+1)/4}
  \\[8mm]
\raisebox{-5mm}{\begin{picture}(10,10)
  \put( 0, 0){\vector( 1, 1){10}}
  \put( 6, 4){\vector( 1,-1){ 4}}
  \put( 4, 6){\line(-1, 1){ 4}}
  \put( 0,10){\makebox(0,0)[br]{$i$}}\put(10,10){\makebox(0,0)[bl]{$j$}}
  \put( 0, 0){\makebox(0,0)[tr]{$k$}}\put(10, 0){\makebox(0,0)[tl]{$l$}}
  \put( 5, 5){\makebox(0,0)[l]{\,\,${\boldsymbol m}$}}
\end{picture}}
\quad &:\,
  (\rarrows_{m}^{-})_{k,l}^{i,j}:=
  \delta_{j,k-m}\delta_{l,i+m}
  \frac{(\qbar)_{i}^{-1}(\qbar)_{j}^{-1}}
       {(\qbar)_{m}(\qbar)_{k}^{-1}(\qbar)_{l}^{-1}}
  (-1)^{i+j+1}
  q^{-ij-(k+l)/2-(N^2+1)/4}
  \\[8mm]
\raisebox{-5mm}{\begin{picture}(10,10)
  \put(10,10){\vector(-1,-1){10}}
  \put( 6, 4){\line( 1,-1){ 4}}
  \put( 4, 6){\vector(-1, 1){ 4}}
  \put( 0,10){\makebox(0,0)[br]{$i$}}\put(10,10){\makebox(0,0)[bl]{$j$}}
  \put( 0, 0){\makebox(0,0)[tr]{$k$}}\put(10, 0){\makebox(0,0)[tl]{$l$}}
  \put( 5, 5){\makebox(0,0)[l]{\,\,${\boldsymbol m}$}}
\end{picture}}
\quad &:\,
  (\larrows_{m}^{-})_{k,l}^{i,j}:=
  \delta_{k,j-m}\delta_{i,l+m}
  \frac{(\qbar)_{i}(\qbar)_{j}}
       {(\qbar)_{m}(\qbar)_{k}(\qbar)_{l}}
  (-1)^{k+l+1}
  q^{-kl-(i+j)/2-(N^2+1)/4}
  \\[8mm]
\raisebox{-5mm}{\begin{picture}(10,10)
  \put( 0,10){\vector( 1,-1){10}}
  \put( 6, 6){\vector( 1, 1){ 4}}
  \put( 4, 4){\line(-1,-1){ 4}}
  \put( 0,10){\makebox(0,0)[br]{$i$}}\put(10,10){\makebox(0,0)[bl]{$j$}}
  \put( 0, 0){\makebox(0,0)[tr]{$k$}}\put(10, 0){\makebox(0,0)[tl]{$l$}}
  \put( 5, 5){\makebox(0,0)[l]{\,\,${\boldsymbol m}$}}
\end{picture}}
\quad &:\,
  (\rarrows_{m}^{+})_{k,l}^{i,j}:=
  \delta_{j,k+m}\delta_{l,i-m}
  \frac{(q)_{i}(q)_{j}}
       {(q)_{m}(q)_{k}(q)_{l}}
  (-1)^{k+l+1}
  q^{kl+(i+j)/2+(N^2+1)/4}
  \\[8mm]
\raisebox{-5mm}{\begin{picture}(10,10)
  \put(10, 0){\vector(-1, 1){10}}
  \put( 6, 6){\line( 1, 1){ 4}}
  \put( 4, 4){\vector(-1,-1){ 4}}
  \put( 0,10){\makebox(0,0)[br]{$i$}}\put(10,10){\makebox(0,0)[bl]{$j$}}
  \put( 0, 0){\makebox(0,0)[tr]{$k$}}\put(10, 0){\makebox(0,0)[tl]{$l$}}
  \put( 5, 5){\makebox(0,0)[l]{\,\,${\boldsymbol m}$}}
\end{picture}}
\quad &:\,
  (\larrows_{m}^{+})_{k,l}^{i,j}:=
  \delta_{k,j+m}\delta_{i,l-m}
  \frac{(q)_{i}^{-1}(q)_{j}^{-1}}
       {(q)_{m}(q)_{k}^{-1}(q)_{l}^{-1}}
  (-1)^{i+j+1}
  q^{ij+(k+l)/2+(N^2+1)/4}
\end{align*}
\par
\medskip
For each local minimum and maximum where an arc labeled with $i$ goes from left to right, we assign the following quantities.
\begin{align*}
\raisebox{-2.5mm}{\begin{picture}(10,5)
  \qbezier(0,5)(1,0)(5,0)
  \qbezier(5,0)(9,0)(10,5)
  \put(6,0){\vector(1,0){0}}
  \put(5,-1){\makebox(0,0)[t]{$i$}}
\end{picture}}
\quad &:\quad
\smile_i:=q^{i-(N-1)/2}
\\[8mm]
\raisebox{-2.5mm}{\begin{picture}(10,5)
  \qbezier(0,0)(1,5)(5,5)
  \qbezier(5,5)(9,5)(10,0)
  \put(6,5){\vector(1,0){0}}
  \put(5,6){\makebox(0,0)[b]{$i$}}
\end{picture}}
\quad &:\quad
\frown_i:=q^{-i+(N-1)/2}
\end{align*}
Then we take the product of all the quantities above and take the summation with all the labels running non-negative integers less than $N$ keeping the labels of two end points of the $(1,1)$-tangle $0$.
(We may choose the labels of the end points arbitrarily.)
\par
Note that $(\darrows_{m}^{\pm})_{k,l}^{i,j}$, $\smile_i$, and $\frown_i$ are obtained from the $R$-matrix $R_J^{\pm}$, $\mu_J$, and $\mu_J^{-1}$ respectively which appear in the enhanced Yang-Baxter operator corresponding to the $N$-dimensional irreducible representation of $sl(2,\C)$
\cite{Kirby/Melvin:INVEM91,Turaev:INVEM88}
(see also \cite{Murakami/Murakami:volume}).
We also note that $\uarrows$, $\larrows$, and $\rarrows$ can be obtained from
$\darrows$, $\smile$, and $\frown$.
\section{The colored Jones function of the figure-eight knot}
Let us consider the figure-eight knot whose $(1,1)$-tangle description is shown in the following figure.
\begin{equation*}
\includegraphics[scale=0.25]{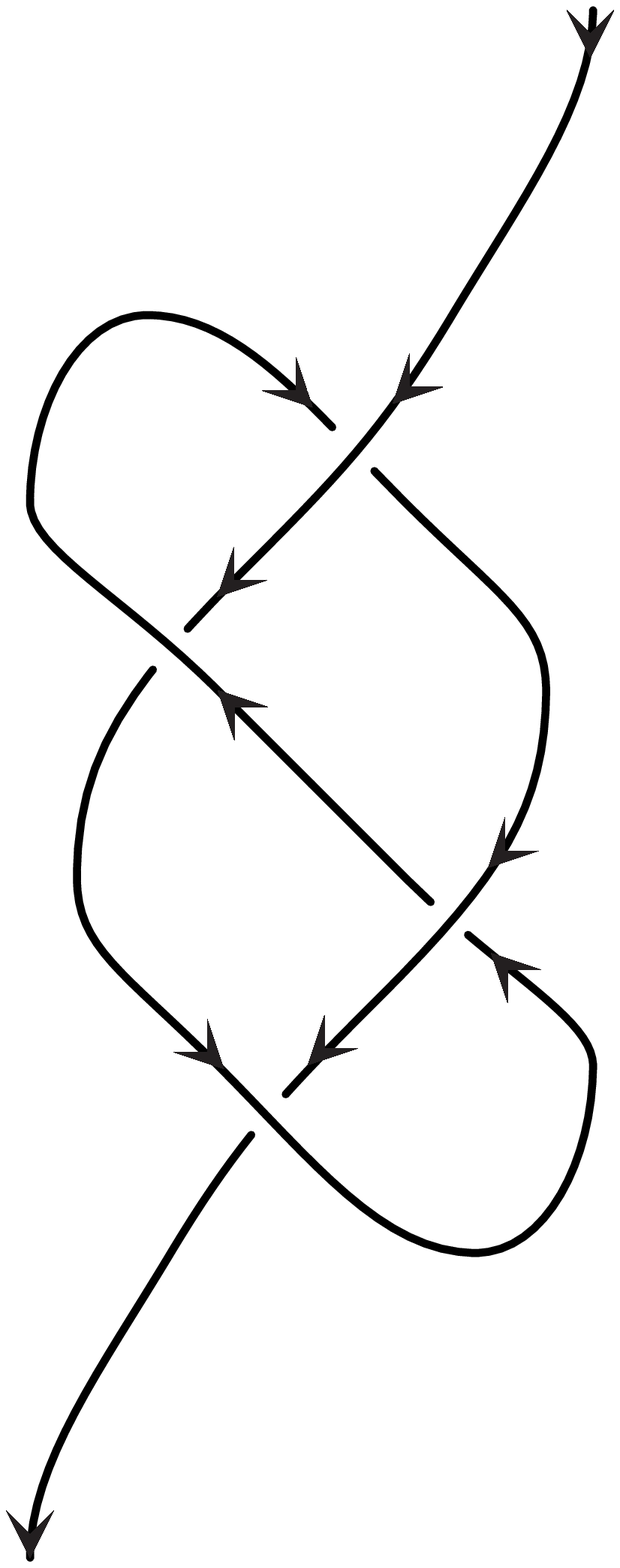}
\end{equation*}
Now we attach a label to each arc and crossing noting that the difference of the labels should be non-negative if we go through an under-crossing and non-positive if we go through an over-crossing
(note Kronecker's delta in $(\darrows_{m}^{\pm})_{k,l}^{i,j}$).
The labels are indicated in the following figure, where integers of bold faces are attached to the crossings.
\begin{equation*}
\includegraphics[scale=0.25]{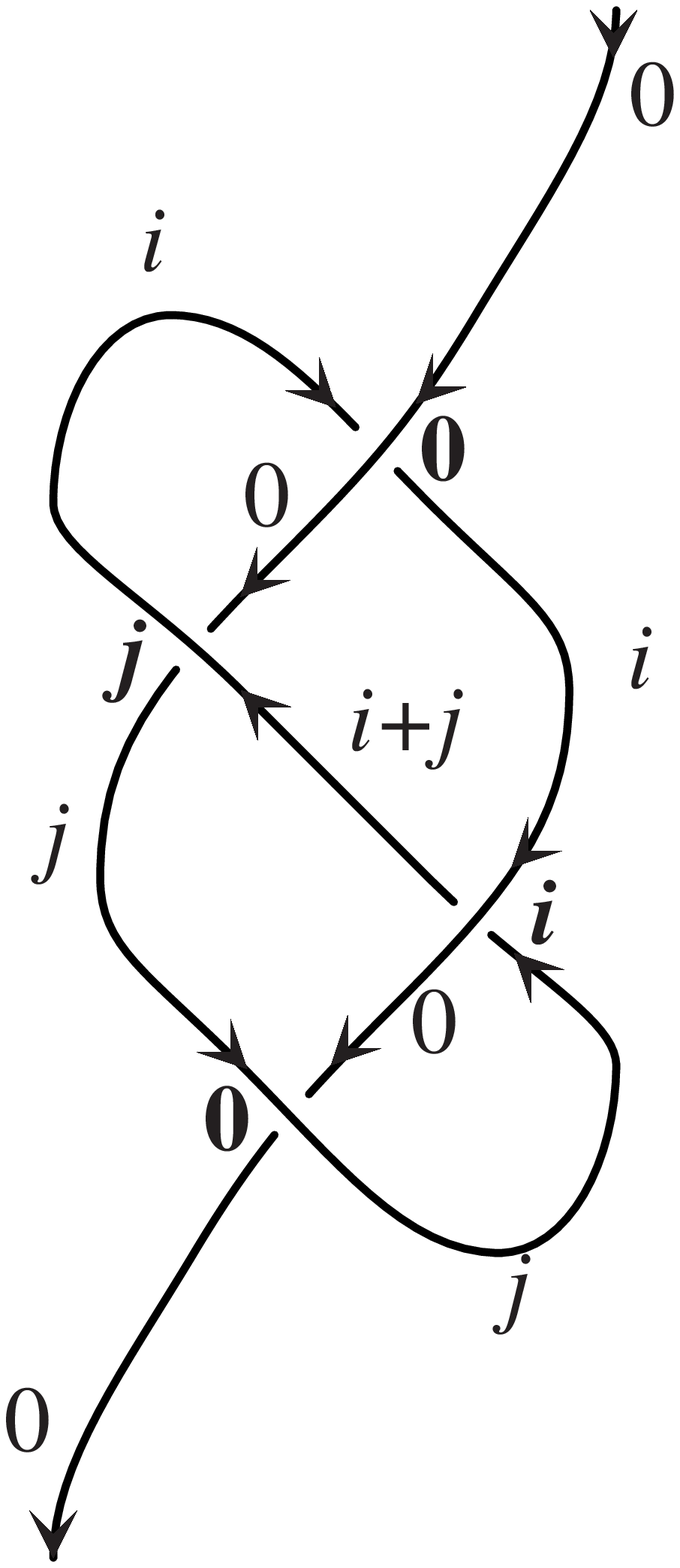}
\end{equation*}
Then our invariant is
\begin{equation}\label{eq:fig8_ij}
\begin{split}
  &J_N(4_1)
  \\
  &\quad=
  \sum_{\substack{0 \le i   \le N-1 \\
                  0 \le j   \le N-1 \\
                  0 \le i+j \le N-1}}
  (\darrows_{0}^{+})_{0,i}^{i,0}\,\times
  (\larrows_{j}^{+})_{j,i+j}^{i,\,\,\,0}\,\times
  (\darrows_{0}^{-})_{0,j}^{j,0}\,\times
  \smile_j\,\times
  (\larrows_{i}^{-})_{\,\,\,0,\,\,\,j}^{i+j,i}\,\times
  \frown_i
  \\
  &\quad=
  \sum_{\substack{0 \le i   \le N-1 \\
                  0 \le j   \le N-1 \\
                  0 \le i+j \le N-1}}
  \frac{(q)_{i}^{-1}}{(q)_{i}^{-1}}(-1)^{i+1}q^{i/2+(N^2+1)/4}
  \frac{(q)_{i}^{-1}}{(q)_{j}(q)_{j}^{-1}(q)_{i+j}^{-1}}(-1)^{i+1}
        q^{i/2+(N^2+1)/4}
  \\
  &\quad\quad\times
  \frac{(\qbar)_{j}}{(\qbar)_{j}}(-1)^{j+1}q^{-j/2-(N^2+1)/4}
  \frac{(\qbar)_{i+j}(\qbar)_{i}}{(\qbar)_{i}(\qbar)_{j}}(-1)^{j+1}
        q^{-j/2-(N^2+1)/4}
  \\
  &\quad\quad\times
  q^{j-(N-1)/2}q^{-i+(N-1)/2}
  \\
  &\quad=
  \sum_{\substack{0 \le i   \le N-1 \\
                  0 \le j   \le N-1 \\
                  0 \le i+j \le N-1}}
  \frac{(q)_{i+j}(\qbar)_{i+j}}{(q)_{i}(\qbar)_{j}}.
\end{split}
\end{equation}
\section{The limit of the colored Jones function of the figure-eight knot}
\label{sec:Ekholm}
In this section I describe the calculation of the limit of
$\log\left(J_N(4_1)\right)/N$ due to T.~Ekholm, and confirm Kashaev's conjecture
for the figure-eight knot.
Putting $k:=i+j$ in \eqref{eq:fig8_ij} we have
\begin{equation*}
\begin{split}
  &J_N(4_1)
  \\
  &\quad=
  \sum_{k=0}^{N-1}\sum_{i=0}^{k}
  \frac{(q)_{k}(\qbar)_{k}}{(q)_{i}(\qbar)_{k-i}}
  \\
  &\quad=
  \sum_{k=0}^{N-1}(-1)^{k}(q^{1/2}-q^{-1/2})^{2k}[k]^{2}
  \\
  &\quad\quad\times
  \sum_{i=0}^{k}
  (-1)^{i}q^{(k^2+k-2ik-2i)/4}\frac{1}{(q^{1/2}-q^{-1/2})^{k}[i]![k-i]!}
  \\
  &\quad=
  \sum_{k=0}^{N-1}(-1)^{k}q^{k(k+1)/4}(q^{1/2}-q^{-1/2})^{k}[k]
  \sum_{i=0}^{k}
  (-1)^{i}q^{-i(k+1)/2}\qbinom{k}{i}
\end{split}
\end{equation*}
since $(q)_l=(-1)^{l}q^{l(l+1)/4}(q^{1/2}-q^{-1/2})^{l}[l]$
and $(\qbar)_l=q^{-l(l+1)/4}(q^{1/2}-q^{-1/2})^{l}[l]$.
Now using \cite[Lemma A.1]{Murakami/Murakami:volume} we have
\begin{equation}\label{eq:fig8_k}
\begin{split}
  J_N(4_1)
  &=
  \sum_{k=0}^{N-1}(-1)^{k}q^{k(k+1)/4}(q^{1/2}-q^{-1/2})^{k}[k]
  (1-q^{-1})(1-q^{-2})\cdots(1-q^{-k})
  \\
  &=
  \sum_{k=0}^{N-1}
  (-1)^{k}
  \left\{(q^{1/2}-q^{-1/2})(q^{2/2}-q^{-2/2})\cdots(q^{k/2}-q^{-k/2})\right\}^2
  \\
  &=
  \sum_{k=0}^{N-1}
  \left\{
    2\sin\left(\frac{1}{N}\pi\right)\times
    2\sin\left(\frac{2}{N}\pi\right)\times\cdots\times
    2\sin\left(\frac{k}{N}\pi\right)
  \right\}^2.
\end{split}
\end{equation}
\begin{rem}
The formula \eqref{eq:fig8_k} was obtained by Kashaev
\cite[(2.2)]{Kashaev:LETMP97}.
The colored Jones polynomial for generic $t$ is
\begin{equation*}
  \sum_{k=0}^{N-1}
  \prod_{l=1}^{k}
  \left(t^{(N+l)/2}-t^{-(N+l)/2}\right)
  \left(t^{(N-l)/2}-t^{-(N-l)/2}\right),
\end{equation*}
which was first obtained by T.~Le.
One can obtain this formula by using technique described here using $R$-matrix for generic $t$.
\end{rem}
Now we will calculate the $N\to\infty$ limit of $\log(J_N(4_1))/N$, which was first obtained by Kashaev \cite{Kashaev:LETMP97}.
The following calculation is due to Ekholm.
\begin{thm}[Kashaev and Ekholm]
\begin{equation*}
  2\pi\lim_{N\to\infty}\frac{\log(J_N(4_1))}{N}=6\Lob\left(\frac{\pi}{3}\right),
\end{equation*}
where $\Lob(\alpha)$ is the Lobachevsky function
\begin{equation*}
  \Lob(\alpha):=-\int_{0}^{\alpha}\log|2\sin\theta|d\theta.
\end{equation*}
Note that the volume of the ideal tetrahedron with face angles $\alpha$, $\beta$, and $\gamma$ $(\alpha+\beta+\gamma=2\pi)$ is
$\Lob(\alpha)+\Lob(\beta)+\Lob(\gamma)$ and that the figure-eight knot complement can be decomposed into two regular ideal tetrahedra.
Therefore the equality above shows that the left hand side equals the volume of
the figure-eight knot, confirming Kashaev's conjecture in this case.
\end{thm}
\begin{proof}
Put $g_k:=\displaystyle\prod_{j=1}^{k}2\sin\left(\frac{j}{N}\pi\right)$ so that
$J_N(4_1)=\displaystyle\sum_{k=0}^{N-1}g_k^2$.
Since
\begin{equation*}
  \begin{cases}
    \displaystyle2\sin\left(\frac{j}{N}\pi\right)
    <1\quad\text{when $\frac{j}{N}<\frac{1}{6}$ or $\frac{j}{N}>\frac{5}{6}$},
    \\[5mm]
    \displaystyle2\sin\left(\frac{j}{N}\pi\right)
    >1\quad\text{when $\frac{1}{6}<\frac{j}{N}<\frac{5}{6}$},
  \end{cases}
\end{equation*}
$g_k$ is decreasing when
$\displaystyle k<\frac{N}{6}$ or
$\displaystyle k>\frac{5N}{6}$ and increasing when
$\displaystyle \frac{N}{6}<k<\frac{5N}{6}$.
Thus (roughly speaking) $g_k$ attains its maximum at
$\displaystyle k=\frac{5N}{6}$.
\par
Since there are $N$ terms in the summation formula of $J_N(4_1)$, we have
\begin{equation*}
  g_{5N/6}^2\le J_N(4_1) \le N g_{5N/6}^2.
\end{equation*}
Taking $\log$ and dividing by $N$ we have
\begin{equation*}
  \frac{2\log g_{5N/6}}{N}
  \le
  \frac{\log J_N(4_1)}{N}
  \le
  \frac{2\log g_{5N/6}}{N}+\frac{\log N}{N},
\end{equation*}
which turns out to be
\begin{equation*}
  2\sum_{j=1}^{5N/6}\frac{\log 2\sin\left(\frac{j}{N}\pi\right)}{N}
  \le
  \frac{\log J_N(4_1)}{N}
  \le
  2\sum_{j=1}^{5N/6}\frac{\log 2\sin\left(\frac{j}{N}\pi\right)}{N}
  +\frac{\log N}{N}.
\end{equation*}
Since $\displaystyle\lim_{N\to\infty}\frac{\log N}{N}=0$, we have
\begin{align*}
    \lim_{N\to\infty}\frac{\log J_N(4_1)}{N}
    &=
    2\lim_{N\to\infty}
    \sum_{j=1}^{5N/6}\frac{\log 2\sin\left(\frac{j}{N}\pi\right)}{N}
    \\
    &=
    2\int_{0}^{5\pi/6}\frac{1}{\pi}\log 2\sin x\,dx
    \\
    &=
    -\frac{2}{\pi}\Lob\left(\frac{5\pi}{6}\right).
\end{align*}
On the other hand from \cite[Lemma 1]{Milnor:BULAM382} we have
\begin{align*}
  \Lob\left(\frac{\pi}{3}\right)
  &=
  \Lob\left(2\frac{\pi}{6}\right)
  \\
  &=
  2\Lob\left(\frac{\pi}{6}\right)+2\Lob\left(\frac{\pi}{6}+\frac{\pi}{2}\right)
  \\
  &\quad\text{(since $\Lob$ has period $\pi$)}
  \\
  &=
  2\Lob\left(\frac{\pi}{6}\right)+2\Lob\left(\frac{\pi}{6}-\frac{\pi}{2}\right)
  \\
  &\quad\text{(since $\Lob$ is an odd function)}
  \\
  &=
  2\Lob\left(\frac{\pi}{6}\right)-2\Lob\left(\frac{\pi}{3}\right).
\end{align*}
Therefore
$\displaystyle\Lob\left(\frac{\pi}{6}\right)
=\frac{3}{2}\Lob\left(\frac{\pi}{3}\right)$
and
\begin{equation*}
  \Lob\left(\frac{5\pi}{6}\right)
  =
  -\Lob\left(\frac{\pi}{6}\right)
  =
  -\frac{3}{2}\Lob\left(\frac{\pi}{3}\right).
\end{equation*}
We finally have
\begin{equation*}
  \lim_{N\to\infty}\frac{\log J_N(4_1)}{N}
  =
  \frac{1}{2\pi}6\Lob\left(\frac{\pi}{3}\right),
\end{equation*}
completing the proof.
\end{proof}
\section{Saddle point method}\label{sec:Kashaev}
In this section I follow Kashaev
\cite{Kashaev:LETMP97}
and calculate the limit of
$\log\left(J_N(4_1)\right)/N$ directly from \eqref{eq:fig8_ij} using the saddle point method.
Note that Kashaev calculated the limit from \eqref{eq:fig8_k} using the same method.
\par
From \cite{Kashaev:LETMP97},
$(q)_i=S_{\gamma}(\gamma-\pi)/S_{\gamma}(-\pi+\gamma+2i\gamma)$
and
$(\overline{q})_i=S_{\gamma}(\pi-\gamma-2i\gamma)/S_{\gamma}(\pi-\gamma)$
with $\gamma=\pi/N$ and $S_{\gamma}(p)$ an analytic function of $p$ which behaves like
\begin{equation*}
  S_{\gamma}(p)\sim
  \exp\left[
        \frac{\Li_2(-e^{p\sqrt{-1}})}{2\gamma\sqrt{-1}}
      \right]
\end{equation*}
for small $\gamma$, where $\Li_2(\zeta)$ is Euler's dilogarithm
\begin{equation*}
  \Li_2(\zeta)=-\int_{0}^{\zeta}\frac{\log(1-\xi)}{\xi}\,d\xi.
\end{equation*}
Putting $z:=q^i$ we have
\begin{align*}
  (q)_{i}
  &=
  \frac{S_{\gamma}(\gamma-\pi)}{S_{\gamma}(-\pi+\gamma+2i\gamma)}
  \\
  &\sim
  \exp
  \left[
    \frac{\Li_2\left(e^{\gamma\sqrt{-1}}\right)
         -\Li_2\left(z\,e^{\gamma\sqrt{-1}}\right)}
         {2\gamma\sqrt{-1}}
  \right]
  \sim
  \exp
  \left[
    \frac{-\Li_2(z)}{2\gamma\sqrt{-1}}
  \right]
\end{align*}
for small $\gamma$.
Similarly putting $w:=q^j$ we have
\begin{align*}
  (\overline{q})_j&\sim\exp\left[\frac{\Li_2(z^{-1})}{2\gamma\sqrt{-1}}\right],
  \\
  (q)_{i+j}&\sim\exp\left[\frac{-\Li_2(zw)}{2\gamma\sqrt{-1}}\right],
  \\
  (\overline{q})_{i+j}&
    \sim\exp\left[\frac{\Li_2(z^{-1}w^{-1})}{2\gamma\sqrt{-1}}\right].
\end{align*}
Therefore for large $N$ we can replace the sum \eqref{eq:fig8_ij} with the following double contour integral
\begin{equation}\label{eq:fig8_double_int}
  \iint
    \exp
    \left[
      \frac{-\Li_2(zw)+\Li_2(z^{-1}w^{-1})+\Li_2(z)-\Li_2(w^{-1})}
           {2\gamma\sqrt{-1}}
    \right]
  dz\,dw
\end{equation}
with suitably chosen contour.
Then the saddle point method (or the method of steepest descent) tells us that
\eqref{eq:fig8_double_int} behaves as a function of $\gamma$ like
\begin{equation}\label{eq:saddle}
  \exp
    \left[
      \frac{F(z_0,w_0)}{2\gamma\sqrt{-1}}
    \right],
\end{equation}
where $F(z,w):=-\Li_2(zw)+\Li_2(z^{-1}w^{-1})+\Li_2(z)-\Li_2(w^{-1})$
and $(z_0,w_0)$ is a solution to the following system of partial differential equations:
\begin{equation}\label{eq:differential}
  \begin{cases}
    \displaystyle\frac{\partial\,F}{\partial\,z}=0,
    \\[5mm]
    \displaystyle\frac{\partial\,F}{\partial\,w}=0.
  \end{cases}
\end{equation}
\par
Since $d\,\Li_2(z)/d\,z=-\log(1-z)/z$, \eqref{eq:differential} turns out to be
\begin{equation*}
  \begin{cases}
    \log(1-zw)+\log(1-z^{-1}w^{-1})-\log(1-z)=0,
    \\
    \log(1-zw)+\log(1-z^{-1}w^{-1})-\log(1-w^{-1})=0
  \end{cases}
\end{equation*}
and so we have
\begin{equation}\label{eq:zw1}
  \begin{cases}
    (1-zw)(1-z^{-1}w^{-1})=1-z,
    \\
    (1-zw)(1-z^{-1}w^{-1})=1-w^{-1}.
  \end{cases}
\end{equation}
Therefore we want to find a solution to the following system of equations:
\begin{equation}\label{eq:zw2}
  \begin{cases}
    z^2w^2-zw-z^2w+1=0,
    \\
    z^2w^2-zw-z+1=0,
   \end{cases}
\end{equation}
giving the trivial solution $z=w=1$.
\par
To get a non-trivial solution we will regard $z$ and $w$ as elements in $\widehat{\C}=\C\cup\{\infty\}$ and put $u:=zw$.
Then \eqref{eq:zw2} becomes
\begin{equation*}
  \begin{cases}
    u^2-u-zu+1=0,
    \\
    u^2-u-z+1=0.
   \end{cases}
\end{equation*}
Throwing away the trivial solution $u=z=1$, we have
\begin{equation}\label{eq:zu}
  \begin{cases}
    z=0,
    \\
    u^2-u+1=0,
   \end{cases}
\end{equation}
giving non-trivial solutions $(z=0, u=\exp(\pi\sqrt{-1}/3))$ and
$(z=0,\exp(5\pi\sqrt{-1}/3))$ ($w=\infty$).
(I learned this `blow-up' technique from J.~Murakami.)
\par
We denote by $F_0$ the $F(z_0,w_0)$
corresponding to the solution $(z=0, w=\infty, u=\exp(5\pi\sqrt{-1}/3))$.
Then
\begin{align*}
  \Im\left(F_0\right)
  &=
  -\Im\left(\Li_2\left(\exp(5\pi\sqrt{-1}/3)\right)\right)
  +\Im\left(\Li_2\left(\exp(-5\pi\sqrt{-1}/3)\right)\right)
  \\
  &\quad
  +\Im\left(\Li_2(0)\right)-\Im\left(\Li_2(0)\right)
  \\
  &=
  2\Im\left(\Li_2\left(\exp(\pi\sqrt{-1}/3)\right)\right),
\end{align*}
where $\Im$ denotes the imaginary part.
Since
$\Im\left(\Li_2\left(\exp(\theta\sqrt{-1})\right)\right)
=\Lob(\theta)+2\Lob(\pi-\theta/2)$
from \cite[(3.7)]{Kashaev:LETMP97} (see also \cite{Kirillov:dilog}),
$\Im\left(F_0\right)$ equals the volume of the figure-eight knot complement.
\par
From \eqref{eq:saddle}
\begin{equation*}
  \lim_{N\to\infty}\frac{\log|J_N(4_1)|}{N}
  =
  \Re\left(\frac{F(z_0,w_0)}{2\pi\sqrt{-1}}\right)
  =
  \frac{\Im\left(F(z_0,w_0)\right)}{2\pi},
\end{equation*}
giving $\Vol(S^3\setminus{4_1})/2\pi$ again, where $\Re$ is the real part.
\section{A cheating calculation}\label{sec:Thurston}
In this section I follow D.~Thurston \cite[page 5]{D.Thurston:Grenoble} to get the limit by a `formal' calculation.
\par
We put
\begin{equation*}
f(i,j):=\frac{(q)_{i+j}(\qbar)_{i+j}}{(q)_{i}(\qbar)_{j}}
\end{equation*}
so that $J_N(4_1)=\sum_{i,j}f(i,j)$.
Now consider the ratios $f(i,j)/f(i-1,j)$ and $f(i,j)/f(i,j-1)$
to find `(local) maxima/minima' of the function $f$ (I am cheating here!).
To do that we will find a solution to the `partial difference equations'
$f(i,j)/f(i-1,j)=f(i,j)/f(i,j-1)=1$, which might give `(local) maxima/minima'.
\par
Since
\begin{equation*}
  \begin{cases}
    \displaystyle
    \frac{f(i,j)}{f(i-1,j)}
    &=
    \displaystyle
    \frac{(1-q^{i+j})(1-q^{-(i+j)})}{(1-q^{i})},
    \\[5mm]
    \displaystyle
    \frac{f(i,j)}{f(i,j-1)}
    &=
    \displaystyle
    \frac{(1-q^{i+j})(1-q^{-(i+j)})}{(1-q^{-j})},
  \end{cases}
\end{equation*}
we have
\begin{equation}
  \begin{cases}
    (1-zw)(1-z^{-1}w^{-1})=1-z,
    \\
    (1-zw)(1-z^{-1}w^{-1})=1-w^{-1}
  \end{cases}
\end{equation}
putting $z:=q^{i}$ and $w:=q^{j}$.
\par
Now we have the same system of equations as in the previous section
\eqref{eq:zw1}.
We denote by $f_{\text{MAX}}$
\footnote{MAX are Nana, Reina, Mina and Lina. Aquarius!}
the $f_{i,j}$
corresponding to the solution $(z=0, w=\infty, u=\exp(5\pi\sqrt{-1}/3))$
($u=zw$).
(I am also cheating since $|q|=1$ and so $z=q^i$ cannot be $0$!)
Since $(q)_i=(\qbar)_j\sim1$ for large $N$ when $q^i=z=0=w^{-1}=q^{-j}$, we have
\begin{align*}
  f_{\text{MAX}}
  &=
  \left|
    (q^{1/2}-q^{-1/2})(q^{2/2}-q^{-2/2})\times\cdots\times
    (q^{(i+j)/2}-q^{-(i+j)/2})
  \right|^2
  \\
  &=
  \left\{
    2\sin\left(\frac{1}{N}\right)\times
    2\sin\left(\frac{2}{N}\right)\times\cdots\times
    2\sin\left(\frac{5\pi}{6}\right)
  \right\}^2
  \\
  &=
  g_{5N/6}^2,
\end{align*}
which is the same value as before and gives the same limit.
\par
Note that a similar calculation using \eqref{eq:fig8_k} was indicated by Thurston and gives the same result as in \S~\ref{sec:Ekholm}.
\section{Geometry}
As seen in \S\S~\ref{sec:Kashaev} and \ref{sec:Thurston},
if we use the triangulation described in
\cite{D.Thurston:Grenoble}
(see the picture below), the only two tetrahedra corresponding to
$u^{-1}=q^{-(i+j)}=\exp(\pi\sqrt{-1}/3)$ survive after taking the limit.
\vspace{10mm}
\begin{equation*}
\raisebox{19mm}{\begin{picture}(10,10)
  \put(10,10){\vector(-1,-1){10}}
  \put( 6, 4){\vector( 1,-1){ 4}}
  \put( 4, 6){\line(-1, 1){ 4}}
  \put( 0,10){\makebox(0,0)[br]{$i$}}\put(10,10){\makebox(0,0)[bl]{$j$}}
  \put( 0, 0){\makebox(0,0)[tr]{$k$}}\put(10, 0){\makebox(0,0)[tl]{$l$}}
  \put( 5, 5){\makebox(0,0)[l]{\,\,${\boldsymbol m}$}}
  \put(80,32){\makebox(0,0)[bl]{$q^{-j}$}}
  \put(43,32){\makebox(0,0)[br]{$q^{k}$}}
  \put(40,-22){\makebox(0,0)[tr]{$q^{i}$}}
  \put(80,-22){\makebox(0,0)[tl]{$q^{-l}$}}
  \put(100,5){\makebox(0,0)[l]{$q^{m}$}}
\end{picture}}
\raisebox{23mm}{\hspace{10mm}$\longrightarrow$\hspace{10mm}}
\raisebox{-4mm}{\includegraphics[scale=0.25]{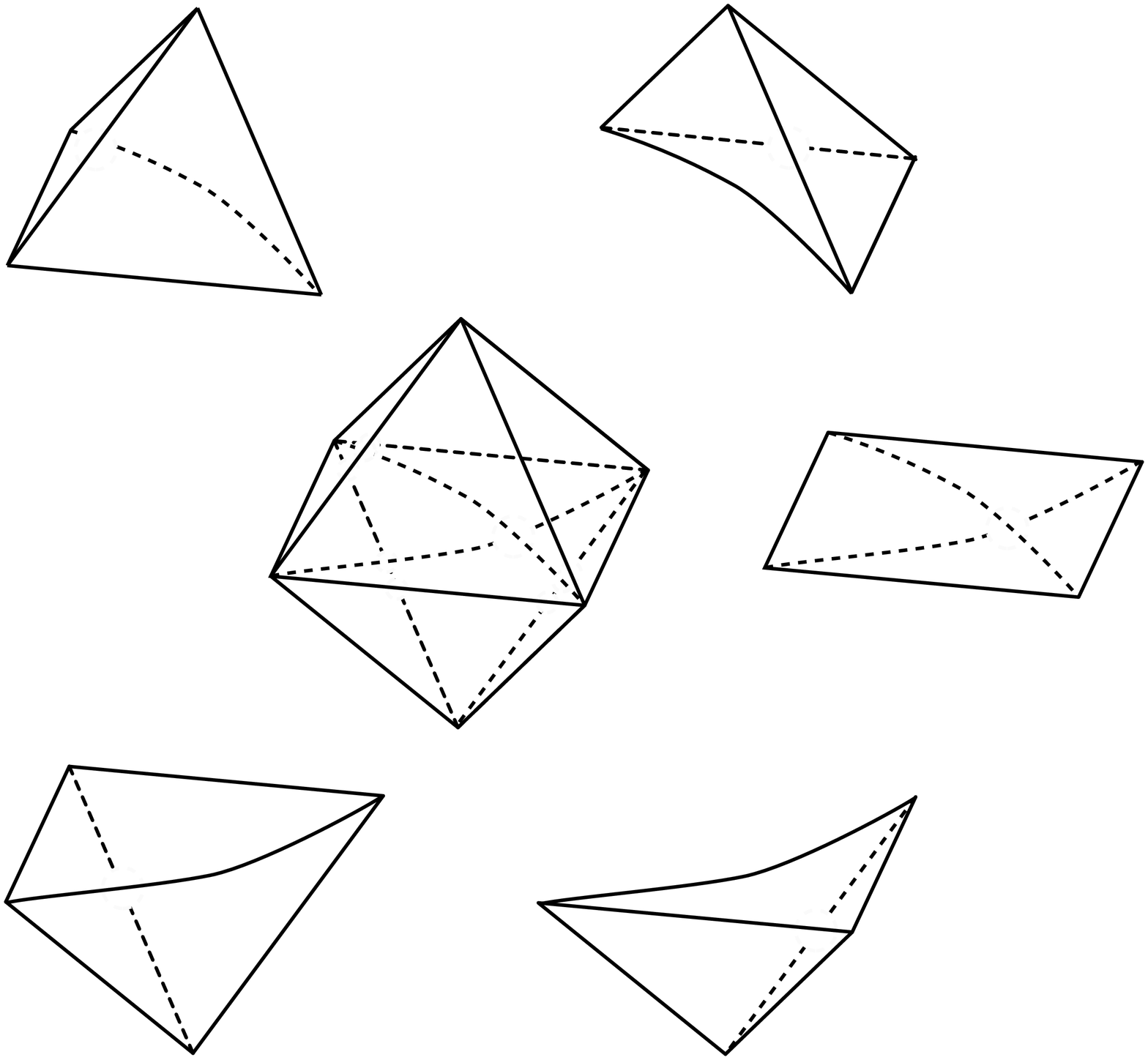}}
\vspace{5mm}
\end{equation*}
Note that the equation $u^2-u+1=0$ in \eqref{eq:zu} is the hyperbolicity equation for the figure-eight knot complement which determines its hyperbolic structure.
\par
The calculations here suggest that in the limit each $R$-matrix corresponds to
five ideal tetrahedra, some of which may collapse, and the partial differential equations (appeared in \S~\ref{sec:Kashaev}) and the partial difference equations (appeared in \S~\ref{sec:Thurston}) give the same algebraic equations \eqref{eq:zu}, which coincide with the hyperbolicity equations.
Due to Yokota this holds in much more general situation
\cite{Yokota:Murasugi70,Yokota:volume}.
Therefore it is now very natural to expect that Kashaev's conjecture is true for any hyperbolic knot and that there should be a rich theory behind it.
\bibliography{mrabbrev}
\bibliographystyle{amsplain}
\end{document}